# Reducing numerical work in non-linear parameter identification

Emoke Imre[1*], Csaba Hegedus[3], Sándor Kovacs[3], Levente Kovacs[2]

[1] HBM Systems research center, Obuda University, Hungary
[2] AIAM Doctoral School, Obuda University, Hungary
[3] Department of Numerical Analysis, Eotvos Lorand University, Budapest, Hungary
[*] Corresponding author, e-mail: imreemok@gmail.com

**Abstract**

The real-life merit functions have an unimaginable complexity of an *M*-dimensional topography, where *M* is the number of the parameters. It is shown that there is an underlying noise-free merit function, called follower merit function which can be constructed from simulated, noise-free data using the solution of a Least Squares minimization. The difference of these is controlled by the norm of the error vector. The local minima arisen from the noise can be skipped during such minimization that apply adjustable, large step sizes, depending on the norm of the error vector.

The suggested hierarchical minimisation - based on the implicit function of a projection map - makes the minimisation in two steps. In the first step some parameters are eliminated with conditional minimisation, resulting in a kind of deepest section of the merit function, called clever section with respect to the remainder parameters. Then the clever section is minimized. The method can be used to handle a quasi-degenerate global minimum, in sensitivity analyses, reliability testing.

## 1 Introduction

### 1.1 Definitions, statements

*1.1.1 Inverse problem*

Let us assume that the model is given by a function $u : R^+ \times R^M \to R$; $(t,\mathbf{p}) \to u(t,\mathbf{p})$; where $R^+ \ni t$ is the time variable, $D \ni \mathbf{p}$ is the parameter vector with $M$ elements; the $D$ is the simply connected and compact parameter domain contained by $R^M$. We can introduce a partial function $v_p(t) = u(t,\mathbf{p})$ and define its graph $G : R^M \to R^2$ as $F = G(\mathbf{p}) = \{R^2 \ni (t, v_p(t))\}$. The function $v_p(t) = u(t,\mathbf{p})$ is called as model solution function.

In the direct problem - for a fixed value of $\mathbf{p}$ - the graph of the partial function $v_p(t) = u(t,\mathbf{p})$, expressed as $F = G(\mathbf{p})$, is determined. In the inverse problem the goal is the determination of the parameter vector $\mathbf{p}$ for a given graph $F$, expressed as $\mathbf{p} = G^{-1}(F)$. In the discrete formulation, only some points of the graph of the partial function $v_p(t)$ are considered, denoted by $v_p(t^i)$, $i=1..N$. These values are compiled into a so-called data vector $\mathbf{f}$.

The measured data vector:

$$\mathbf{f} = \mathbf{u_m} = \begin{bmatrix} f(t^1) \\ f(t^2) \\ .... \\ f(t^N) \end{bmatrix} \quad (1)$$

where $t_i$ are the sampling times ($i=1,2..N$).

Simulated data vector using a fixed value of the parameter vector:

$$\mathbf{f'} = \begin{bmatrix} u(t^1, \mathbf{p_*}) \\ u(t^2, \mathbf{p_*}) \\ .... \\ u(t^N, \mathbf{p_*}) \end{bmatrix} \quad (2)$$

The function $\mathbf{u}(\mathbf{p}): R^M \to R^N$, $(\mathbf{p}) \to \mathbf{u}(\mathbf{p})$, called model response function, containing $N$ elements, can be constructed from the copies of $u(t,\mathbf{p})$ in such a way that times $t_1...t_N$ - the elements of the so-called sampling time vector $\mathbf{t}$ - are substituted for the time variable $t$ in $u(t,\mathbf{p})$.

$$* u(p) = \begin{bmatrix} u(t^1, p) \\ u(t^2, p) \\ \cdots \\ u(t^N, p) \end{bmatrix} \quad (3)$$

In the inverse problem, for a fixed value of *f* and *t*, the goal is the determination of the parameter vector *p*, from the equation $h(p) = f - u(p) = 0$, expressed symbolically as $p = u^{-1}(f)$.

$$h(p) = f - u(p) = \begin{bmatrix} f(t^1) \\ f(t^2) \\ \cdots \\ f(t^N) \end{bmatrix} - \begin{bmatrix} u(t^1, p) \\ u(t^2, p) \\ \cdots \\ u(t^N, p) \end{bmatrix} = 0 \quad (4)$$

The solution of a system of non-linear of equations $f - u(p) = 0$ is needed, which generally has weak solution only, since generally $N > M$. The weak solution is the global minimizing vector $p_{\min}$ of the least squares objective function called real-life merit function which is the norm squares of the error vector $h(t,p) = u_m(t_i) - u(t, p)$:

$$F(p) = h(p)^T h(p) = min! \quad (5)$$

*1.1.2 The general linear inverse problem:*

The general linear inverse problem is to fit a set of data points to a model that is not just a linear combination of 1 and *t* (namely a + b *t*), but rather a linear combination of any M specified functions of *t*. The general form of this kind of model is as follows:

$$u(p) = \begin{bmatrix} f_1(t) & f_2(t) & \cdots & f_k(t) \end{bmatrix} \begin{bmatrix} p_1 \\ p_2 \\ \cdots \\ p_M \end{bmatrix} = \sum_i p_i f_i(t) \quad (6)$$

In the general linear inverse problem, the model response function:

$$u(p) = \begin{bmatrix} f_1(t^1) & f_2(t^1) & f_k(t^1) \\ f_1(t^2) & f_2(t^2) & f_k(t^2) \\ f_1(t^N) & f_2(t^N) & f_k(t^N) \end{bmatrix} \begin{bmatrix} p_1 \\ p_2 \\ p_M \end{bmatrix} = \begin{bmatrix} u(t^1, p) \\ u(t^2, p) \\ \cdots \\ u(t^N, p) \end{bmatrix} \quad (7)$$

The error vector and the direct formulation of the inverse problem:

$$h(p) = f - u(p) = \begin{bmatrix} f(t^1) \\ f(t^2) \\ \cdots \\ f(t^N) \end{bmatrix} - \begin{bmatrix} u(t^1, p) \\ u(t^2, p) \\ \cdots \\ u(t^N, p) \end{bmatrix} = 0 \quad (8)$$

$$\begin{bmatrix} f_1(t^1) & f_2(t^1) & f_k(t^1) \\ f_1(t^2) & f_2(t^2) & f_k(t^2) \\ f_1(t^N) & f_2(t^N) & f_k(t^N) \end{bmatrix} \begin{bmatrix} p_1 \\ p_2 \\ p_M \end{bmatrix} = A\ p = f = \begin{bmatrix} f(t^1) \\ f(t^2) \\ \cdots \\ f(t^N) \end{bmatrix} \quad (9)$$

The Gauss Normal Equations:

$$A^T A\ p = A^T f \quad (10)$$

The merit function:

$$F(p) = (f - u(p))^T (f - u(p)) = \sum_{i=1}^{N} (f(t^i) - u(t^i, p))^2$$
$$= (A\ p - f)^T (A\ p - f) \quad (11)$$

*Statement*

In the general linear inverse problem, the second derivative of the merit function with respect to *p* is globally positive definit, being a constant Gram matrix.

*Proof*

The derivative is a constant Gram matrix (see Eq 9):

$$F(p)'' = 2 A^T A \quad (12)$$

*1.1.3 The partly linear inverse problem*

The parameter space is split into the direct sum of two subspaces:

$$p = [p_1, p_2], \quad p_1 \in D_1 \in R^k, \ p_2 \subset D_2 \subset R^{M-k} \quad (13)$$

The partly linear inverse problem is to fit a set of data points to a model that is a linear combination of any M specified functions of *t* and $p_2$. The general form of this kind of model is

$$u(p) = \begin{bmatrix} f_1(t, p_2) & f_2(t, p_2) & \cdots & f_k(t, p_2) \end{bmatrix} \begin{bmatrix} p_{1,1} \\ p_{1,2} \\ \cdots \\ p_{1,k} \end{bmatrix} = \sum_i p_{1,i} f_i(t, p_2) \quad (14)$$

In the partly linear inverse problem, the model response function:

$$u(p) = \begin{bmatrix} f_1(t^1, p_2) & f_2(t^1, p_2) & f_k(t^1, p_2) \\ f_1(t^2, p_2) & f_2(t^2, p_2) & f_k(t^2, p_2) \\ f_1(t^N, p_2) & f_2(t^N, p_2) & f_k(t^N, p_2) \end{bmatrix} \begin{bmatrix} p_{1,1} \\ p_{1,2} \\ p_{1,k} \end{bmatrix} = A(p_2) p_1$$

(15)

In the partly linear inverse problem, the error vector and the direct formulation of the conditional inverse problem:

$$h(p) = \begin{bmatrix} f(t^1) \\ f(t^2) \\ \cdots \\ f(t^N) \end{bmatrix} - \begin{bmatrix} u(t^1, p) \\ u(t^2, p) \\ \cdots \\ u(t^N, p) \end{bmatrix} = 0 \quad (16)$$

In the partly linear inverse problem, the Gauss Normal Equations of the conditional inverse problem for $p_1$, at fixed $p_2$:

$$A(p_2)^T A(p_2) p_1 = A(p_2)^T f \quad (17)$$

In the partly linear inverse problem, the merit function can be written as follows:

$$F(p) = (f - u(p))^T (f - u(p)) = \sum_{i=1}^{N} (f(t^i) - u(t^i, p))^2 = (A(p_2)p_1 - f)^T (A(p_2)p_1 - f) \quad (18)$$

*Statement*
In the partly linear inverse problem, the second derivative of the merit function with respect to $p_1$ is globally positive definit, being a constant Gram matrix at the fixed value of $p_2$.

*Proof*
The derivative with respect to $p_1$ is a constant Gram matrix in case of the fixed value of $p_2$ (see Eq 15):.

$$F(p)_{p_1}'' = 2 A(p_2)^T A(p_2) \quad (19)$$

**1.2 Difficulties in non-linear parameter identification**

**1.2.1 Problems in non-linear minimisation**

The multidimensional, non-linear minimization genetally needs iterations and it is generally difficult if the model has several unknown parameters since (i) too many local minima of the real-life merit function may occur due to the noise, moreover, (ii) the geometry of the valley of the global minimum could be asymmetric and (iii) the minimum can be quasi-degenerated.

The first problem - the unimaginable complexity of an *N*-dimensional topography [1] - is treated in the literature by some ways. A trial and error process is used in every case, when not significant decrease in the error is attained, the process may restart from a new starting vector.

The methods, which include one-dimensional minimisation iteration steps, are generally coupled with some new starting vector determination methods (eg,. the simulated annealing method, which steps away from a local minimum using probability theory methods, see eg. [1] or by fitting a smooths hypersurface locally [2]).

The non-linear minimization methods without one-dimensional minimisation (eg., simplex method [1], or secant method see Section 6) are considered as effective tool in case of a few parameters [1], the reason of the better performance is related to the fact that these may skip some critical points, as shown here.

The quasi-degenaret minimum problem is treated generally by regularization ([3]), the asymmetric shape of the vally is an unsolved problem.

**1.2.2. Problems in reliability testing**

The singularity of the Hesse matrix at a minimum minimum $p_{\min}$ gives information on local uniqueness. If the global minimum is degenerated, infinite many solutions may exists, some parameterts are dependent. In case of linear inverse problems, the global uniqueness is met in case of local non-singularity.

The solution of a real-life inverse problem $p_{\min}$ is approximate with some error (see the Appendix). Using the linear model, the computation of the standard deviation of the parameters is possible, assuming normally distributed noise (residuals) [1]. The geometrical interpretation is related to the orthogonal projection of a confidence domain, which is defined by a level line of a kind of noise-free merit function (Fig. 1), with symmetry related to the individual parameters. The unsolved problem is the parameter error in case of non-normal noise.

In case of non-linear inverse problems, the global uniqueness is not met in case of local non-singularity. The computation of the standard deviation is possible, similarly to the linear case, assuming normally distributed noise [1]. It has similar geometrical meaning: it is the orthogonal projection of the confidence domain, related to a noise-free, linearized merit function. Some unsolved, additional questions are the global uniqueness, and the parameter error in case of non-normally distributed noise and the non-symmetry with respect to the individual parameters.

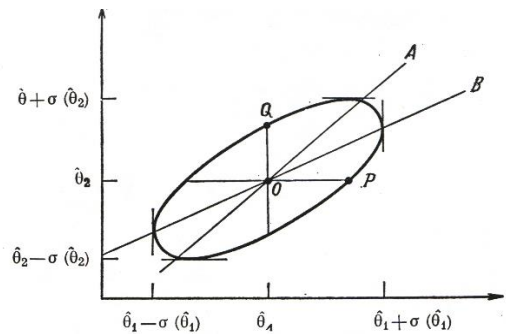

**Fig. 1** Confidence domain of 68.3%, the standard deviation in case of normally distributed noise (after [1]).

**1.3 The content of paper**

To overcome the foregoing difficulties, three new concepts are presented and discussed. The first concept is the Hierarchical Inverse Problem, which means split minimisation, such that different minimisation methods may be used in each part.

The hierarchical inverse problem concept is mathematically explained by the implicit function theorem. The related map is the orthogonal projection of the graph of the merit function onto a plane $<F, x>$, and the implicit function $y = g(x)$ is related to the fold of the image of the map.

The second concept is related to the use of the section of the merit function along the fold, defined in the HIP concept. The so-called clever section can be used to estimate the global minimum without minimization if some independent information for the $x_{min}$ part of the global minimum is available. In this case $y_{min} = g(x_{min})$. The one-dimensional clever section - depending on a single parameter - is a deepest sensitivity section which can be used for sensitivity analyses and reliability testing.

The third concept is the follower, noise-free merit function - the best possible noiseless approximation of the real-life merit function - which can equally be used for minimisation and reliability testing. The real-life merit function has too many local minima with nearly the same merit function value around the global minimiser $p_{min}$. The global uniqueness of the solution can be interpreted by the follower here. The obstruction domain for minimisation is a parameter error domain, defined by a level line of the follower merit function here.

The follower is not known a priori but can be determined afterwards using the global minimiser $p_{min}$. It is shown that the difference of the two functions is dependent on the actual size of the noise-free residuals. Due to the similarity, the decreasing direction of the follower can be estimated during minimisation by the noisy (real-life) merit function on a part of the parameter domain - aiding minimisation - if the step size is large enough, eg. dependent on the actual size of the noise-free residuals.

In the paper two examples are presented. The step size of the secant method, a non-linear minimization methods, without one-dimensional minimisation ([6 to 8]) is analysed. It is shown that this method is theoretically able to skip some local minima since its step size is dependent on the actual size of the noise-polluted residuals. The example also illustrates how the linearly dependent parameters can be eliminated as an application of the Hierarchical minimisation concept.

The second example shows an automatic global minimisation related to the injective solutions of linear PDE-s, giving reliability information as well. The the injective solution was expressed with 8 parameters, 4 non-linearly dependent parameters („model solution H") and, in its variants the number of the non-linearly dependent parameters was decreased up to 1 by using information concerning the analytical properties of the model (see eg., „model solution HC"). The number of the function value evaluation depended on the number of non-linearly dependent parameters exponentially, in case of a single non-linearly dependent parameters, it was less than 60, increasing by 1 order of magnitude with an additional non-linearly dependent parameter. The case of 3 non-linearly dependent parameters was numerically comparable with the use of the conjugate gradient method coupled with a new starting vector determination method, but only the former was informative on the reliability of the solution of the inverse problem.

**3 Hierarchical inverse problem solution**

**3.1 General**

*3.1.1. Definition*

The parameter space can be split into the direct sum of two subspaces:

$$p = [p_1, p_2], \quad p_1 \in D_1 \subset R^J, p_2 \subset D_2 \subset R^{M-J} \quad (20)$$

The first inverse problem is a conditional minimization:

$$F(p_1, p_2) = \min!, \quad p_2 = \text{const}, \forall p_2 \in D_2 \subset R^{M-J} \quad (21)$$

the solution of which is the relation $p_1 = a(p_2)$. If this inserted then an M-J dimensional minimal section of $p_2$ the section of the original Least Squares merit function, denoted by $F[a(p_2),(p_2)]$ is resulted which depends only on $p_2$.

The second inverse problem:

$$F[a(p_2), p_2] = \min! \quad p_2 \in D_2 \subset R^{M-J} \quad (22)$$

is the minimization of $F(a(p_2),(p_2))$.

*3.1.2 Definition*

The section $F[a(p_2),(p_2)]$ is called minimal or clever section of $p_2$.

### 3.1.3 Statement

The solution of the hierarchical minimisation is trivially equivalent to the original one if the Hessian of the Least Squares merit function $F$ with respect to $p_1$ is globally positive definite, which is met e.g., if the function $F$ is strictly convex, or in case of linear dependence.

## 3.2 Mathematical explanation

We change the notation so that subscripts can be omitted in the index. We use $y$ i $x$ instead of $p_1$ and $p_2$.

### 3.2.1 The definition of the projection map

Let us consider a strictly convex, non-negative, analytic Least Squares merit function $F(p)$: $D \subset \mathbf{R}^{n+m} \to \mathbf{R}^1$ with simply connected, compact domain $D$ and with the arbitrary $p=(x, y)$ arbitrary direct sum decomposition of $D$. The global minimum is $(a,b)$, where $F$ has the only critical point.

Let us consider the orthogonal projection of the graph of $F$ $\{(x, y), F(x, y)\}$ onto the $F$-$x$ coordinate plan in the total Euclidean space generated by the graph. The projection planes are $x=$**const** coordinate planes, the critical values of the map - where the fold of the projection is found - are related to the $F_y$' $(x, y) = f(x, y) = \mathbf{0}$.

### 3.2.2 The implicit function theorem in case of the a prejction of the graph a strictly convex function

The implicit function theorem [4, 5] is applied in the case of the projection of the graph of a strictly convex, non-negative, analytic (Least Squares merit) function to describe the fold where the condition is met $F_y = f(x, y) = \mathbf{0}$ as follows. Starting from a direct product decomposition, the goal is to construct a function $g$: $\mathbf{R}^n \to \mathbf{R}^m$ whose graph ($\mathbf{x}$, $g(\mathbf{x})$) is precisely the set of all ($\mathbf{x}$, $\mathbf{y}$) such that the $F_y' = f(x, y) = \mathbf{0}$.

Let us consider the critical point at the global minimum $(a,b) = (a_1, ..., a_n, b_1, ..., b_m)$ where $F_y'$ $(x, y) = f(x, y) = \mathbf{0}$, the zero is element □ of $\mathbf{R}^m$. If the matrix $F_y$ $_y$'' $(x, y) = [(\partial f_i/\partial y_j)(a, b)]$ has full rank, then there exists an open set $U$ containing $a$, an open set $V$ containing $b$, and a unique, analytic, injective function $g$: $U \to V$ such that

$$\{\mathbf{x}, \mathbf{g}(\mathbf{x}) | \mathbf{x} \in U\} = \{(\mathbf{x},\mathbf{y}) | \mathbf{x} \in U, \mathbf{f}(\mathbf{x},\mathbf{y}) = \mathbf{0}\} \quad (23)$$

It follows from the injectivity, that if parameter vector $x_{\min}$ is the part global minimum, then $y_{\min} = g(x_{\min})$.

### 3.2.3 Consequences

Consequences of the strict convexity are as follows.

The implicit function is globally unique if $F_{yy}$'' $(x, y) = [(\partial f_i/\partial y_j)(a, b)]$ is globally positive definite, which is met due to the strict convexity.

The implicit functions related to various $\mathbf{x'} < \mathbf{x} < \mathbf{p}$ have the same lattice structure like the lattice structure of subspaces of the $R^M$ generated by the direct product decomposition related to $\mathbf{x'} < \mathbf{x} < \mathbf{p}$. If the split is such that $\mathbf{x'} < \mathbf{x}$ then the clever section $F[g(\mathbf{x}), \mathbf{x}]$ contains $F[g(\mathbf{x'}), \mathbf{x'}]$. The $F[g(\mathbf{x'}), \mathbf{x'}]$ can equally be determined from $F[g(\mathbf{x}), \mathbf{x}]$ or $F(\mathbf{p})$.

## 4 Follower merit function

### 4.1 General

#### 4.1.1 The definition of the follower merit function

The so-called follower merit function $F(p)$ is defined as the norm squares of the simulated error vector, computed by using the simulated, noise-free data, generated with $p_{\min}$:

$$F'(\mathbf{p}) = \sum_{i=1}^{N} h_i(\mathbf{p})^2 = \sum_{i=1}^{N} [u(t_i, \mathbf{p}) - u(t_i, \mathbf{p}_{min})]^2 \quad (24)$$

where $u(t_i, p_{\min})$, the simulated, noise-free data and $u(t, p)$ is the model solution, $t_i$ ($i=1..N$) are sampling times, the notation is not related to derivation, it means noise-free. The dash is used to indicate simulated noise-free data.

#### 4.1.2 The definition of the noise

The noise $z_i$ can be expressed from the following decomposition, after the determination of the global minimum of the merit function and the simulated data:

$$u_m(t_i) = u(t_i, \mathbf{p}_{min}) + z_i \quad (25)$$

where $p_{\min}$ minimiser, $u(t, p)$ is the model solution.

#### 4.1.3 The definition of the similarity domain and the parameter error domain

The level surface of the noise-free merit function related to the global minimum of the real-life merit function is used to define the parameter error domain, similarly to the statistical error domain (Fig. 1). The parameter error domain is defined as the follows, by a sublevel set of the follower merit function:

$$E = \{p \in D | F'(p) \leq c = F(p_{min})\} \quad (26)$$

In the parameter error domain the real-life merit function contains several local minima with about the same value. In the complement domain - far from the global minimum - is the following:

$$\|h'(p)\|^2 > \|z\|^2 \quad (27)$$

is met. This domain is defined as the so-called similarity domain, where the follower merit function can be approximated by the real-life merit function.

### 4.2 The relation of the two merit functions

*Statement*

Using the foregoings, the two merit functions: the follower merit function $F$ and the real-life merit function $F$ has the following relationship:

$$F(p) = \sum_{i=1}^{N}[h_i(p) - z_i]^2 = F'(p) - \sum_{i=1}^{N}[2h_i(p)z_i] + F(p_{min}) \quad (28)$$

By using the Cauchy-Schwartz-Bunyakovszkij inequality:

$$F'(p) - 2\|h'(p)\| \|z\| + F(p_{min}) = \|h'(p)\|^2 - 2\|h'(p)\| \|z\| + \|z\|^2$$
$$\leq F(p) \leq \|h'(p)\|^2 + 2\|h'(p)\| \|z\| + \|z\|^2 = \quad (29)$$
$$F'(p) + 2\|h'(p)\| \|z\| + F(p_{min})$$

The difference of the two functions can be expressed as follows:

$$-2\|h'(p)\| \|z\| + \|z\|^2 \leq F(p) - F'(p) \leq 2\|h'(p)\| \|z\| + \|z\|^2 \quad (30)$$

According to Eq (30), the difference of the two function seems to be small for fixed, small noise, being controlled however by the norm of the noise-free error vector.

## 5. Practical use of concepts and methods

### 5.1 Use of split minimisation

*5.1.1 Elimination of the linearly dependent parameters*
The most important practical application of the split minimisation is the elimination of the linear part of parameter vector through sub-minimisation, being the related Hesse matrix globally positive definite. The linear subminimisation - a generalized linear inverse problem [1] - can be solved without iteration, eg., by the SVD or singular value decomposition (SVD) algorithm. The linearly dependent parameters can be eliminated by sub-minimisation before each function value evaluation in case of any non-linear minimization method.

*5.1.2   Use of the clever sections*

The implicit function $g(x)$ can be used to estimate the global minimiser of $F[x,g(x)]$ without minimization in case the value of $x_{min}$ is available from independent information. In this case $y_{min} = g(x_{min})$. This follows from the fact that the implicit function is injective.

The one-dimensional clever sections of the follower merit function and the real-life merit functions can simultaneously be represented in the 2-dimensional plane, reflecting uniqueness and parameter error (Fig.2). Being the error domain non-symmetric, the orthogonal projection of its diameter onto parameter axis $p^i$ is reasonable to be used to characterize the parameter error.

### 5.2 Minimisation with the follower
*5.2.1 General*
Statement
The follower and the real-life merit functions are pseudo-metrics. The metrics related to the real-life merit function $F$ is introduced as follows:

$$d_F(a,b) = |F(a) - F(b)| \quad (31)$$

where $D \ni a,b$ are parameter vectors. In case $b$ is the global minimum, the merit function indicates the distance from the global minimum, being a pseudo-metrics.

Proof
The definition of a metric on the set $\{\alpha, \beta, \gamma....\}$ is as follows:

1)  $d(\alpha,\beta) \geq 0$,  2)  $d(\alpha,\beta) = 0 \Leftrightarrow \alpha = \beta$,  3)  $d(\alpha,\beta) = d(\beta,\alpha)$,  4)
$$(32)$$

If the second is not met, a psudometric is defined.

*5.2.2   Suggested methods*
Consequence
According to Eq (30), the difference between the follower and the real-life merit functions is bounded by a term which is linearly related to the norm of the noise-free, follower at a fixed noise vector. If the norm of the noise is relatively small, then theoretically (i) the follower merit function values can be approximated by real-life merit function values, (ii) decreasing directions can be assessed, (iii) the step size can be linked to the difference between the follower and the real-life merit functions.

Step sizes
In the non-linear parameter identification methods not using line minimization, the step size can theoretically be large enough to skip the critical points due to the noise. This condition is observed in Example 1.

Coordinate grid

In case of strictly convex follower merit function, some decreasing convex level surfaces (and a minimum) may approximately be bracketed on the similarity domain, by coordinate hyperplanes related to the points of a global coordinate grid, using real-life merit function values.

The same global coordinate grid can be used - with a projection algorithm (Fig. 3)) - to determine the one-dimensional clever sections. Once the minimum of the real-life merit function is determined, the same algorithm can be used to determine the one-dimensional clever sections of the follower merit functions, too. These can be used to determine the error domain (Fig. 2) and, to observe the difference of the follower and real-life merit functions in the similarity domain.

A local coordinate grid can be generated in the subspace of some important parameters only to find a decreasing direction and a new starting vector in this direction.

The grid dimension is determined by the number of the non-linearly dependent parameters, since the linearly dependent parameters can be eliminated. The real-life merit functions values are computed in the points of a coordinate grid in a nested cycle. The one-dimensional clever section point of each parameter $x_i$ can be determined by projections along each respective coordinate planes $x^{i,j}$ =constant, such that for a given parameter value $x^{i,j}$ all grid values are considered and, the parameter values related to the conditional minimum is selected (Fig. 3).

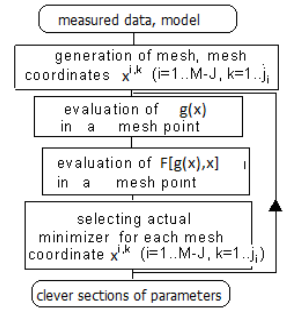

Fig. 3. The determination of the approxinmate one-dimensional clever sections by elimination, projection and minimum search in a nested cycle system. A coordinate mesh is generated in the M-J dimensional space of the nonlinearly dependent parameters, J linearly dependent parameters are eliminated.

## 6 Example 1 - secant method

In the secant-type methods the step size is controlled by the real-life merit function $F$ value (i.e., norm square of the noise-polluted error vector), and a similar variable basically controls the difference of the two functions, also. These methods are shortly presented as follows. The notation is slightly changed here, the parameter number is $n$ and the number of sampling times is $m$.

### 6.1 Wolfe method

In the secant method [6, 7, 8] $n+1$ trial solutions for the parameter vector and the related residual or error vector ($p_i$, $h(p_i)$) ($i = 0,...,n$) are used to compute a new parameter vector:

$$\begin{bmatrix} \boldsymbol{p} \\ \boldsymbol{h(p)} \end{bmatrix} \approx \begin{bmatrix} \boldsymbol{p_0} \\ \boldsymbol{h(p_0)} \end{bmatrix} + \begin{bmatrix} \boldsymbol{A}_\mathrm{p} \\ \boldsymbol{A}_\mathrm{r} \end{bmatrix} \boldsymbol{q} \qquad (33)$$

where $m$ is the number of sampling times, $n$ is parameter number. $A_\mathrm{p}=[p_0- p_1, ..., p_0- p_n]$, $A_\mathrm{r}=[h(p_0)- h(p_1)-, ..., h(p_0)- h(p_n)]$ and $q^\mathrm{T}=[q^1, ..., q^n]$ is a vector with $n$ coefficients of secant method, called as an auxiliary variable.

The $q$ in the one-dimensional case is equal to $h(p_1)/[h(p_0) - h(p_1)]$, similar is valid in multidimensional case for tzhe elements of $q$. The absolute value of the multiplier $h(p_1)$ is a pseudometric, expressing a distance from the global minimum. The step size is comparable with the controlling term of difference between the follower and the real-life merit function (see Eq (10)). Far from the global minimum, a large step from $p_0$ is resulted, allowing the skipping of some critical points. In other words, the in the secant method, the step size (related to the $q$ coefficient) depends on the same variable as the difference between the follower and the real-life merit functions. These norms or merit function values are pseudometrics, being large far from the minimum and vice versa.

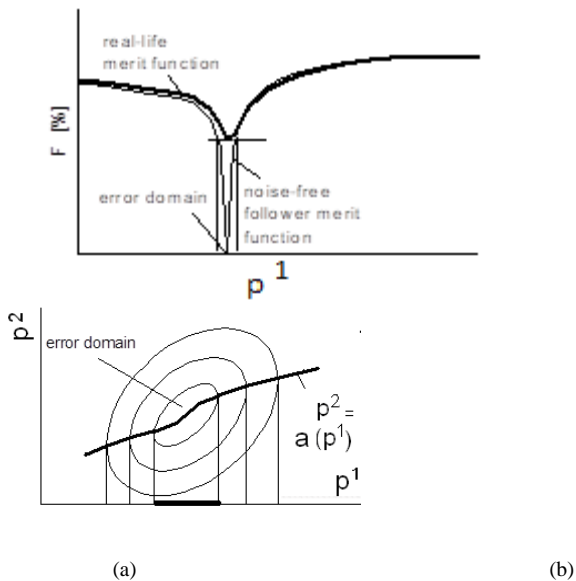

(a)          (b)

Fig. 2. The one dimensional clever section of merit functions and the parameter error domain (a) in the total space (b) in the parameter space. Being the error domain non-symmetric, the orthogonal projection of its diameter onto parameter axis $p^i$ is reasonable to be used. See its value in App.

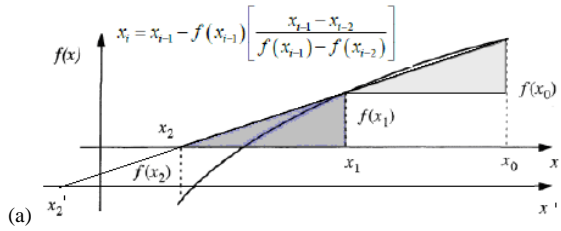

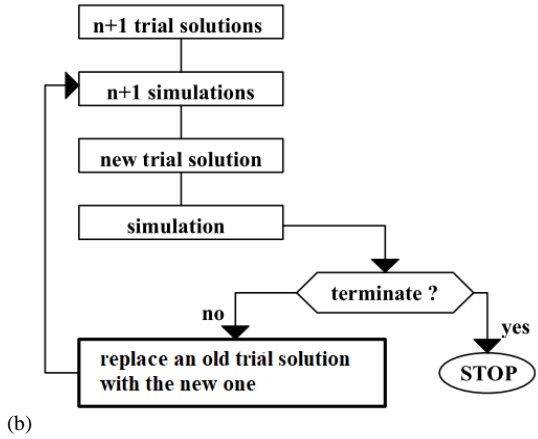

**Fig. 4.** (a) Secant method in 1D. Note that if the function value at *x*1 is the larger, then the step x2-x1 is the greater (see by the shift of the x axis). (b) Flow diagram of Wolfe-secant method

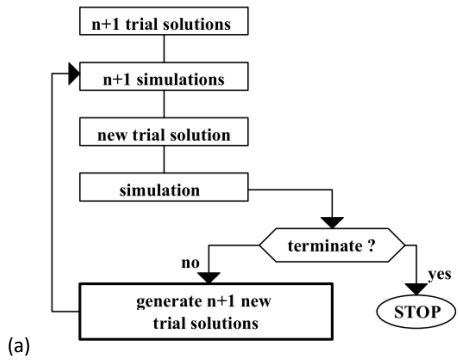

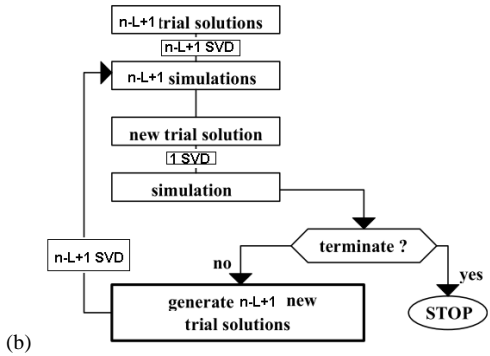

**Fig. 5**. Flow diagrams (a) Modified secant method (b) The completion of the modified shark secant method with linear subminimisation (each function value evaluation is preceeded by a linear sub-minimisation). The question is that the linear subminimisation may amend the method being the secant method linear in each iteration cycle.

*6.2 The modified Wolfe method*
Wolfe ([5]) suggested to select $n+1$ new trial solutions for the next iteration by dropping the trial solution for which the norm of the error vector is the largest. This procedure fails when the norm of the new residual vector is the largest one. In the amended method ([2], Fig. 5), after determining $\boldsymbol{p}_{n+1}$, a second trial solution $\boldsymbol{p}_{n+2}$ is determined. With matrix equation:

$$\begin{bmatrix} T_P & 0 \\ 0 & T_r \end{bmatrix} \begin{bmatrix} A_p \\ A_r \end{bmatrix} q_2 = \begin{bmatrix} p_0 - p_{n+1} \\ h(p_0) \end{bmatrix} \qquad (34)$$

where $T_p$ and $T_r$ are transformation diagonal matrices with elements:

$$t_p^i = \frac{p_{n+1}^i - p_{n+2}^i}{p_0^i - p_{n+1}^i}, t_r^i = h^i(\boldsymbol{p}_{n+1})/h^i(\boldsymbol{p}_0) \qquad (35)$$

The distance of $\boldsymbol{p}_{n+1}$ and $\boldsymbol{p}_{n+2}$ is related to the decrease of the error and the merit function between $\boldsymbol{p}_{n+1}$ and $\boldsymbol{p}_0$. The $n+1$ new trial solutions are selected from the new vectors as:

$$\boldsymbol{p}^i = \boldsymbol{p}_{n+1} + \boldsymbol{e}_i \left( \boldsymbol{p}_{n+2}^i - \boldsymbol{p}_{n+1}^i \right) \qquad (36)$$

where $\boldsymbol{e}^i$ is the i-th unit vector. The system of linear equations in the second rows of the matrix equations are solved by singular value decomposition (SVD).

The new auxiliary variable $q_2$, is determined from the first line of the matrix equation (17). In the one-dimensional case, it is a multiple of $q$, being equal to $h(p_1)/h(p_2)*q$ [8]. A similar generalization is valid in the multidirmnsional case, describen in [8]. It can be noted that the absolute value of the multiplier is a pseudometric again, expressing a distance from the global minimum, and the norm of this term is found in Eq ( 30). Far from the global minimum, it is large, a large step from $\boldsymbol{p}_1$ is determined.

*6.3 Elimination of the linearly dependent parameters*
The most important practical application of the HIP in minimisation is the elimination of the linear part of parameter vector through sub-minimisation. The linear subminimisation can be solved by singular value decomposition (SVD) which can be added to any algorithm. If the linear part of the parameter vector is eliminated in The modified secant method (Fig. 5) then the numerical work is smaller since the number of the parameters is *n-L* instead of *n* where *L* is the number of the linearly dependent parameters. The size of matrices $A_p$ and $A_l$ is smaller. However, before each function value evaluation a linear sub-minimisation using the SVD algorithm is needed.

**Table 1** Types and notation of one dimensional oedometric dissipation tests with constant boundary condition

| (Multistage) oedometric relaxation test | MRT or ORT |
|---|---|
| (Multistage) oedometric compression test | MCT or OCT |

**Table 2** Types dissipation tests made with static penetrometers, modelled with cylindrical and spherical (ellipsoid) shaped domain.

| Measured variable | dissipation test made by |
|---|---|
| Pore water pressure dissipation test, sensor on the shaft and/or on the tip | CPTu (static cone penetrometer) |
| Total stress dissipation test, sensor on the shaft and/or on the tip | CPT- piezo-lateral stress cell, CPT- DMT tip, CPT - CPT- qc |
| Effective stress dissipation test, sensor on the shaft | CPT- fs |

Table 3 Point-symmetric consolidation models [12]

| $V$ or $\varepsilon$ boundary condition | 1D (Oedometric models) | 2D (Cylindrical pile models) | 3D (Spherical pile models) |
|---|---|---|---|
| no (uncoupled) | Terzaghi (1923) | Soderberg (1962) | Torstensson (1975) |
| $v$-$v$ (coupled 1) | Imre (1997-1999) | Imre & Rózsa (1998) | Imre & Rózsa (2002) |
| $v$-$\varepsilon$ (coupled 2) | Biot (1941) | Randolph at al (1979) [ | Imre & Rózsa (2005) |

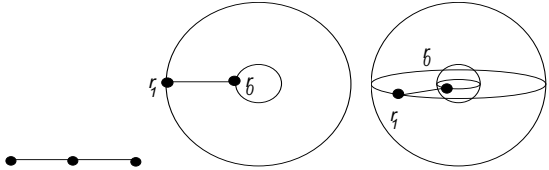

**Fig. 6.** The displacement domain bounded by a (a) 0 dimensional sphere (oedometer model), (b) 1 dimensional sphere (cylindrical model), (c) 2 dimensional sphere (spherical model).

## 7. Example 2 - dissipation test evaluation

The objective of the study was to elaborate automatic evaluaton methods for such dissipation tests of Geotechnics that can be evaluated with linear, point-symmetric, coupled consolidation models [9 to 16] (Tables 1 to 3).

The follower merit function related to these injective, analytical consolidation model solutions, with fixed, convex initial condition and sufficiently long data series is strictly convex, according to the experiences [11].

However, by the identification of the initial condition on a function class allowing non-convex analytic functions and by modelling of relaxation (time dependent constitutive law), this geometry may slightly change, but remains near-convex. The same method - applying the global coordinate grid projection algorithm - was applicable in every case (software is available for each case [17 to 22]).

In case of the standard multistage oedometric compression tests (MCT) and its dual counterpart, the multistage oedometric relaxation tests (MRT), a short testing procedure (with shorter stages than 99% of the dissipation time except for the last stage) was suggested. Laboratory tests were made on identical soil samples with different plasticity [9, 10]. The model validation showed that (i) the fitting error decreased to half for those model-variants where the simultaneous creep/relaxation was considered, (ii) the relaxation qualitatively influenced the model response in a parameter dependent way, and (iii) the short stage data were successfully evaluated using the clever section of the single parameter $c$.

In this paper a summary is given on the evaluation of the multistage oedometric relaxation tests MRT. It is shown, that by using the suggested technique, a class of inverse problems - being related to (near-)convex follower merit functions - becomes easily soluble without iteration or regularization, the reliability criteria of the solution can be readily be tested.

### 7.1 Model-versions

*7.1.1 Analytical solution*

The analytical solution of the linear, point-symmetric coupled consolidation model of the oedometric relaxation test for the so-called total stress $\sigma$ and the so-called pore water pressure $u$ (Appendix A to C [8 to 11]):

$$\sigma(t) = -\sum_{i=1}^{\infty} \varphi_i e^{-i^2 \pi^2 T} + \sigma_\infty \qquad (37)$$

$$u(t,y) = -\sum_{i=1}^{\infty} \varphi_i [\cos(y(i\pi/H)) - 1] e^{-i^2 \pi^2 T} \qquad (38)$$

where $T = ct/H^2$ is the time factor, $c$ is the coefficient of consolidation, $\phi_i$ are Fourier coefficients related to the initial condition, $H$ is the half-width of a double-drained sample (length between the boundaries of the model).

The linear, coupled consolidation model consolidation model has some special features, which can be summarized as follows ([8 to 11]). Being the volume constant during a stage, the mean effective stress ($\sigma_{mean}$) is equal to the final total stress ($\sigma_\infty$) which depends on the displacement load $v_0$ and the $H$ the sample height:

$$\sigma_\infty = \frac{E_{oed} v_0}{H} \equiv \sigma'_{mean}(t) \qquad (39)$$

The time dependent constitutive law was included in an approximate form into the model in its most complex form, by adding a relaxation term to the total stress. The relaxation term includes a partial unloading effect. The equations used for model fitting were as follows.

$$\sigma(t) = -\sum_{i=1}^{\infty} \varphi_i e^{-i^2 \pi^2 T} + \sigma_\infty - \Delta\sigma^r(t) \qquad (40)$$

The empirical relaxation term:

$$\Delta\sigma^r(t) = s\sigma(0)\frac{1}{1-sb}\log\left(\frac{t+t_1}{t_1+t_3}\right); t > t_3 \qquad (41)$$

where $s$ is the coefficient of relaxation, $t_3$ is the pause of relaxation in case of partial unloading, $b=\log((t_1+t_3)/t_1)$, $t_1$ is time parameter.

*7.1.2 Model-versions for parameter identification*
The initial condition was identified using the following parametric shape function in model H:

$$u(0,y) = A\ y^3 + B\ y^2 + C\ y \qquad (42)$$

where $A$, $B$ and $C$ are real valued initial condition parameters, with linear dependence. The related Fourier coefficients were included in the model.

The model H contained 8 parameters, 4 out them were with non-linear dependence (Table 4). The linearly dependent $A$, $B$ and $C$ parameters equally allowed convex, concave or partly convex and concave initial condition functions. Several model-versions were elaborated, differing in the number of the non-linearly dependent parameters, using additional information on the consolidation part-model and the relaxation part-model (described in Appendix B and C), in here only model-version HC is mentioned where no relaxation was modeled.

It can be noted that parameters $\sigma_\infty$ and $c$ have different meaning if relaxation is modelled or if not. If the relaxation is modelled, the compression curve constructed from the identified $\sigma_\infty$ is time independent, this the so-called instant compression curve. However, if the relaxation is not modelled (ie., for model-version HC), the identified $\sigma_\infty$ is depending on the stage duration, as it is dependent in the measurements (as it is well-known from the oedometer testing).

**7.2 Algorithms for parameter identification**

*7.2.1 Conjugate gradients with a new starting vector algorithm*

The model-version H was fitted on measured data with the conjugate gradient method. However, the solution was not possible to be determined due to the large number of critical points. Therefore, a local grid system was used to determine some new starting vectors. The grids were 2-dimensional, considering the most important parameters ($c$, initial condition). The conjugate gradient method was used with 5 to 8 parameters.

*7.2.2. Projection algorithm*

The model-version H and HC were fitted on measured data with the projection algorithm including linear subminimisation [8 to 11], in the frame of an automatic global minimization (Fig. 3). The global coordinate grid was generated in the subspace of the non-linearly dependent parameters, the linearly dependent parameters were eliminated by subminimisation. The grid size was defined by using the half error domain diameter values of some preliminary simulations, regarding the clever section of noisy and the follower, noise-free merit functions.

The one-dimensional clever sections were approximately determined by a projection algorithm, using real-life merit function values in the points of a coordinate grid. Once the minimum of the real-life merit function was determined, the same algorithm was reused to determine the one-dimensional clever sections of the follower merit functions, too. These were used to determine the error domain and, to observe uniqueness under the effect of added non-linearities.

**7.3 Results**

*7.3.1 New starting vector method*

Three successive iteration steps, including the conjugate gradient starting and ending points and coarse grids centered the ending critical points determined by the conjugate gradient method are shown in Figure 7. The section along the fold of the projection within the local grids were represented, revealing the complex geometry of the merit function (Fig. 7). The success of the local grid was apparent but no reliability information was found for the solution and the computational work was excessive.

### 7.3.2. Spit inverse problem solution (projection algorithm and clever sections)

The projection algorithm was used with linear subminimisation according to Fig. 3. In case of model-version H a 4-dimensional coordinate grid was used, the four one-dimensional clever sections were determined by projection and a subsequent minimization for of the real-life and the follower merit functions. In case of model-version HC, 1-dimensional grid was used to determine the clever section of the real-life and the follower merit functions for parameter $c$.

The clever clever section of parameter $c$ determined for the two model-versions were similar, indicating two minima both models (Fig. 8), one minimum was related to a nice, physically admissible initial condition, another to a physically non-admissible one.

The compression curve points ($\sigma_\infty$) were determined by model HC from each stage. In case of short stages, the good minimum was quasi-degenerated (Fig. 9). Evaluating the short stages with the simplest model HC, the implicit function $g(c)$ was determined for each stage only. From function $g(c)$, by using the $c$ identified from the last, long stage, the global minimiser vector and parameters $\sigma_\infty$ were selected. The latter gave compression curve points, being situated above the one determined by the long compression test in accordance with the expectations following from the time dependency of the constitutive law (Fig. 9).

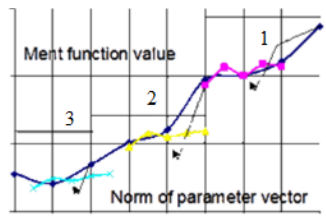

**Fig. 7.** The complexitiy of the topography, implicit functions at three successive steps, fine and coarse local grids. Conjugate gradient method (with arrows) ([16]).

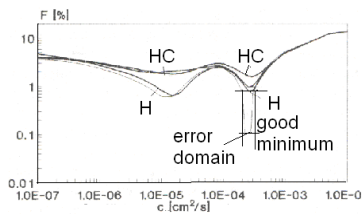

**Fig. 8.** Model H and HC, long stage, clever sections c. Note the non-unique solution for both models ([14]).

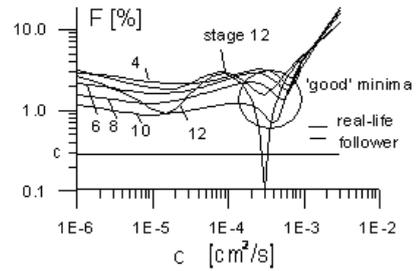

**Fig. 9.** Model HC, short stages 1 to 11 and last, long stage 12. The clever sections $c$ is quasi-degenerated for stages 1 to 11. The c identified from the long stage 12. was used for the evaluation of the short stages ([14]).

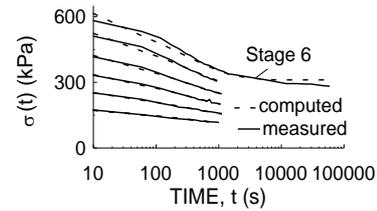

(a) Total stress.

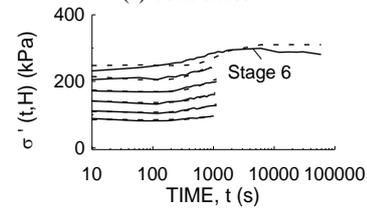

(b) Effective stress.

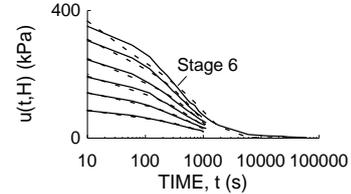

(c) Pore water pressure

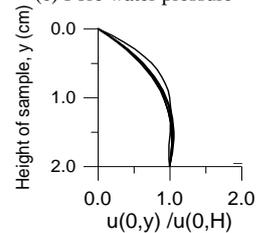

(d) Identified initial condition in the various stage.

**Fig.10.** MRT, fitted and measured stresses, identified initial condition in the various stages.

**Table 4.** Parameters of the most general oedometric model-version H

| Model | Parameter | Symbol | Dependence |
|---|---|---|---|
| Consolidation | Initial condition | A, B , C | linear |
| | Coeff. of consolidation | $c$ | nonlinear |
| | asymptotic total stress | $\sigma_\infty$ | linear |
| Relaxation | coefficient of relaxation | $shH$ | nonlinear |
| | pause of relaxation | $t_3$ | nonlinear |
| | delay time | $t_1$ | nonlinear |

Table 5. Parameters of the consolidation oedometric relaxation test model-version HC

| Model | Parameter | Symbol | Dependence |
|---|---|---|---|
| Consolidation | initial condition | A, B, C | linear |
| | Coeff. of consolidation | $c$ | nonlinear |
| | asymptotic total stress | $\sigma_\infty$ | linear |

Table 6. Number of function value evaluations, oedometric relaxation test, model-version HC and H

| Model-version | H | HC |
|---|---|---|
| Parameter # (linear) | 4 | 4 |
| Parameter # (non-linear) | 4 | 1 |
| # of global mesh points | 125000 | 54 |

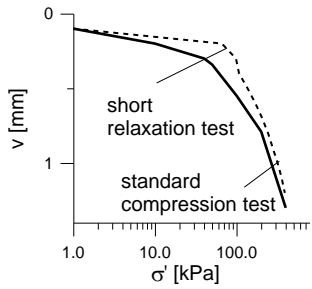

**Fig. 11.** Compression curves measured by standard compression tests and identified from short stage MRT data, using model HC. The $c$ identified from the last stage was used for the evaluation of the short stages with quasi-degenerate minimum (see Fig. 9).

## 8. Discussion, conclusion

The results can be summarized and discussed as follows.

### 8.1 The hierarchical minimisation of analytic Least Squares merit functions

The non-linear parameter identification related to a strictly convex, analytic Least Squares merit function can be split in terms of dimension. In the first part of the split inverse problem, a section of the merit function ($F[g(\mathbf{x}), \mathbf{x}]$, called clever section ) is determined by subminimisation, then the clever section is minimised. Two different minimisation methods can be used at two parts of the direct sum decomposition of the parameter domain. The method can be applied repetitively.

Geometrically, the fold of such a map is determined which is the orthogonal projection of the graph of the merit function onto a coordinate plane $<F, \mathbf{x}>$, where $p=(\mathbf{x},\mathbf{y})$ is a direct sum decomposition. The implicit function $g(\mathbf{x})$ of the derivative of the projection with respect to $\mathbf{y}$ can be determined by sub-minimization at every fixed $\mathbf{x}$.

Analytically, the multi-dimensional implicit function theorem was applied for the $\mathbf{y}$-derivative of a projection map related the graph of a merit function onto a vertical coordinate plane, generated by the complement $\mathbf{x}$, where $p=(\mathbf{x},\mathbf{y})$ is a direct sum decomposition of the parameter space. The implicit functions - for strictly convex merit functions - are globally unique and have the same lattice structure as the subspaces of the direct sum decomposition.

The method can be used for not only for strictly convex, analytic Least Squares merit function but also in the genral case on condition that the Hessian with respect to $\mathbf{y}$ is globally positive definite. The most important practical application is the elimination of the linear part of parameter vector through sub-minimisation, which can be added to any algorithm. The linearly dependent parameters can be eliminated by sub-minimisation before each function value evaluation.

### 8.1 The clever section of the analytic Least Squares merit functions

The one-dimensional clever section $F[g(x_i), x_i]$ of the merit function related a parameter $x_i$ can be used to determine the sensitivity of the parameter and the diameter of the error domain, moreover, to test uniqueness. The one-dimensional clever sections can approximately be determined by a projection algorithm suggested here.

An important practical application could be the case of the quasi-degenerate minimum. The $g(x)$ can be used to estimate the global minimum of $F[x,g(x)]$ without minimization it, on condition that independent knowledge is available for $x_{min}$. By using the known parameter vector value $x_{min}$ of the global minimum, since $y_{min} = g(x_{min})$. This follows from the fact that the implicit function is injective.

### 8.3 Follower merit function

The real-life merit functions have an unimaginable complexity of an $M$-dimensional topography, where $M$ is the number of the parameters. However, there is an underlying noise-free merit function, called follower merit function which can be constructed from simulated, noise-free data using the solution of a Least Squares minimization.

The follower merit function is not known a priori but can be constructed once the global minimiser is known. In this way it can be used in reliability testing, can be used to define a parameter error domain and to define uniqueness criteria. The follower merit function may have much less critical points, can be even convex.

Concerning convexity, the following comment can be made. There is a false opinion that no bracketing procedure can be used in multidimensional space [1].

However, there is an ideal exception which is the class of the is strictly convex minimisation. Hyperplanes generated by coordinate a mesh points can bracket decreasing convex level surfaces in multidimensional case.

### 8.4 Examples

The relatively small difference between the follower, noise-free and the real-life, noisy merit functions is controlled by the value of the follower noise-free merit function, which is a pseudo-metric, (expressing the distance from the global minimum, where its value is zero). In the secant method, the step size is linearly dependent on the merit function value. There are several consequences. The step size is comparable with the size of the irregularities of the real-life merit function with respect to the follower merit function on the similarity domain. If the function value at $x_1$ is large, far from the global minimum, the stepsize is large and vice versa; the global minimum is approached effectively. In addition, being the method inherently linear, it may occur that a linear model part is treated without the inclusion of a separate, linear algorithm for the elimination.

In the case of strictly convex, follower merit functions, the bracketing of the level surfaces is theoretically possible on the $M$-dimensional topography by the hyperplanes of a properly selected coordinate grid (generated in the parameter space). In case of real-life inverse problems, only the follower, noise-free merit function can be strictly convex due to the noise and, the convexity can be used on the similarity domain.

According to the experiences, this may be the case in case of the injective solutions of linear system of PDE-s, in case of fixed, convex initial condition and sufficiently long data series. In case of the injective solutions of linear, coupled point-symmetric consolidation models (Table 3), the follower merit function can be strictly convex. Some non-iterative, non-linear parameter identification methods were elaborated related to these injective solutions, by using the suggested projection algorithm and global coordinate grid. In the second example, an important practical application was illustrated on the treatment of the quasi-degenerate minimum in case of "short data series".

### 8.5 further research

Further research is suggested, among others
- on the convexity of the follower merit functions related to injective solutions of system of coupled, linear partial differential equations; including the effect of non-convex initial conditions and an added non-linearity;
- on the difference between real-life and convex follower merit functions, by using one-dimensional clever sections;
- on the differences between the various (statistical and geometrical) parameter error domains in the light of the true error sizes;
- on such step size of the non-linear, multidimensiona minimization methods that allows the use the follower merit function in minimisation;
- on the application of the hierarchical minimisation of Least Squares merit function in case of system of non-linear equations and in the case of general non-linear optimization where the merit functions may have a slightly different type, for example it can not necessarily be expressed in the form of a Least Squares merit function defined in section 1;
- on the secant – type non-linear parameter identification methods, the performance in case of linear models parts of general non-linear merit function and Least Squares merit functions.

**Appendix**

Example for the introduced definitions at Least Squares inverse problem

Let us assume the following model with a single linearly dependent parameter $a$:

$$u(t,x) = at^2 x \tag{43}$$

where t is the time variable $x$ is the space variable (fixed in here), $a$ is the single parameter.

Simulated data vector assuming $x=1, t=[1,2,3]$, $a = 2$ (i.e. fixed values for the sampling times, the space variable and parameter $a$):

$$f' = \begin{bmatrix} 2 \\ 8 \\ 18 \end{bmatrix} \tag{44}$$

Noisy data vector and noise vector:

$$f_m = \begin{bmatrix} 2.1 \\ 7.8 \\ 18.2 \end{bmatrix}, \tag{45}$$

$$z = \begin{bmatrix} 0.1, \\ -0.2, \\ 0.2 \end{bmatrix},$$

The direct formulation of the noisy inverse problem:
$$Aa = f_m \tag{46}$$

The direct formulation of the noise-free problem:
$$Aa = f \tag{47}$$

where matrix $A$ is as follows (with elements being equal to the derivative of the model with respect to $a$ at fixed value of the sampling time vector and the space variable):

$$A = \begin{bmatrix} 1 \\ 4 \\ 9 \end{bmatrix} \tag{48}$$

In detail, the direct formulation of the noise-free inverse problem:

$$\begin{bmatrix} 1 \\ 4 \\ 9 \end{bmatrix} a = \begin{bmatrix} 2 \\ 8 \\ 18 \end{bmatrix} \tag{49}$$

In detail, the direct formulation of the noisy inverse problem

$$\begin{bmatrix} 1 \\ 4 \\ 9 \end{bmatrix} a = \begin{bmatrix} 2.1 \\ 7.8 \\ 18.2 \end{bmatrix} \tag{50}$$

The solution can be determined by the Gauss Normal Equations:

$$A^T A \quad a' = A^T f_m \tag{51}$$

Solution with the Gauss Normal Equations for the noise-free problem:

$$A^T A \quad a = A^T f \tag{52}$$

In detail:
$$A^T A = 98 \quad A^T f = 196 \quad a = 2 \tag{53}$$

Solution with the Gauss Normal Equation for the noisy problem:

$$A^T A = 98 \quad A^T f_m = 197.1, \quad a' = 2{,}011224 \qquad (54)$$

Simulated data vector assuming $x=1, t=[1,2,3], a=2.011224$ (i.e. fixed values for the sampling times, the space variable and parameter $a$):

$$f' = \begin{bmatrix} 2.011224 \\ 8.044896 \\ 18.101016 \end{bmatrix} \qquad (55)$$

Real-life merit function:

$$F(a) = (2.1-a)^2 + (7.8-4a)^2 + (18.2-9a)^2 \qquad (56)$$

The simulated, noise-free merit function related to then true value of the parameter $a$:

$$F'(a) = (2-a)^2 + (8-4a)^2 + (18-9a)^2 \qquad (57)$$

The follower, simulated noise-free merit function, related to the solution of the noisy inverse problem, the identified value of parameter $a$:

$$F'(a) = (2{,}011224-a)^2 + (8.044896-4a)^2 + (18.101016-9a)^2 \qquad (58)$$

The clever sections of the noisy and the follower merit function for parameter a are shown in Fig. 2. The true error (~0,011) is equal to about 1/5 of the width of the error domain (~0,056) or about 4/10 of the half width of the error domain called generalized standard deviation.

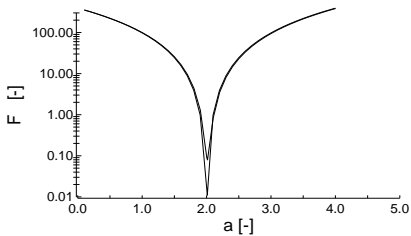

**Fig.** 12. The clever sections of the noisy and the follower merit function for parameter *a*

**Appendix A The consolidation model**
The point-symmetric, linear coupled consolidation models can be summarized into a single mathematical model in the function of the embedding space dimension (Fig 6, Table 3, oedometric m=1, cylindrical m=2, and spherical models m=3). For every embedding space there are two models differing in one boundary condition. In this section the linear coupled consolidation models elaborated for the relaxation/compression test models (m=1) are presented, denoted by ORT and OCT models here.

*1) System of differential equations*
Two equations have been derived from the equilibrium condition and, from the continuity condition, assuming weightless soil and water. The Equilibrium Equation compiles the equilibrium condition, the effective stress equality, the geometrical and, the constitutive equations, as follows:

$$E_{oed} \frac{\partial \varepsilon}{\partial y} - \frac{\partial u}{\partial y} = 0 \qquad (59)$$

and, the Continuity Equation compiles the continuity equation, the Darcys law and the geometrical equation, as follows:

$$-\frac{k}{\gamma_v} \Delta u + \frac{\partial \varepsilon}{\partial t} = 0 \qquad (60)$$

where $v$ is displacement, $u$ is pore water pressure (neglecting the gravitational component of the hydraulic head), $y$ and $t$ are space and time co-ordinates respectively, $E_{oed}$ is the oedometric modulus. Tthe volumetric strain and the Laplacian operator, for dimension $m=1$ are as follows:

$$\varepsilon = \frac{\partial v}{\partial y}, \Delta = \frac{\partial}{\partial y^2} \qquad (61)$$

*2) Boundary conditions*
The displacement domain of the oedometric test is bounded by two zero dimensional circles, with radii $r0 = 0$ and $r0 = H$, *respectively. For a double-drained sample, the upper-lower surface coincides with the outer boundary, respectively, the symmetry point coincides with the inner boundary. In the following four boundary conditions are presented. Three are common, one is different for the two models.

(1) The (common) *boundary condition Nr. 1*:

$$u(t,r)|_{y=0} = 0 \qquad (62)$$

(2) The (common) *boundary condition Nr. 2*:

$$\frac{\partial u(t,y)}{\partial y}\bigg|_{y=H} \equiv 0 \qquad (63)$$

(3) The (common) *boundary condition Nr. 3*:

$$v(t,y)|_{y=H} \equiv 0 \qquad (64)$$

*(4) Boundary condition Nr. 4* concerning the coupled 1 models:

$$v(t,y)|_{y=0} \equiv v_0 > 0 \qquad (65)$$

(5) *Boundary condition Nr. 5* concerning the compression test models:

$$\varepsilon(t,r)|_{y=0} \equiv \varepsilon_0 > 0. \qquad (66)$$

*3) Analytical solution*

Steady-state solution part
The solution for the relaxation test model:

$$v^p(y) = v_0 \left(1 - \frac{y}{H}\right) \qquad (67)$$

$$\sigma^p(y) = -\frac{E_{oed} v_0}{H} \quad (68)$$

and, for the compression test model:

$$v^p(y) = \frac{\sigma_0}{E_{oed}}(H - y) \quad (69)$$

$$\sigma^p(y) = \sigma_0 \quad (70)$$

Transient solution part

The displacement $vt$ for the oedometric relaxation test model:

$$v^t(t,y) = \sum_{k=1}^{\infty} a_k \cdot \sin\left(\frac{k \cdot \pi}{H} y\right) \cdot e^{-\frac{k^2 \pi^2}{H^2} c_v \cdot t} \quad (71)$$

where $a_k$ ($k=1...\infty$) are the Fourier coefficients of an odd initial displacement function and, for the compression test model:

$$v^t(t,y) = \sum_{k=1}^{\infty} b_k \cdot \cos\left(\frac{(2k-1) \cdot \pi}{2H} y\right) \cdot e^{-\frac{(2k-1)^2 \pi^2}{4H^2} c_v \cdot t} \quad (72)$$

where $b_k$ ($k=1...\infty$) are the Fourier coefficients of an even initial displacement function. The function $u$ for the relaxation test model:

$$u(t,y) = \sum_{k=1}^{\infty} \alpha_k \left\{\left[\cos\left(\frac{k \cdot \pi}{H} y\right)\right] - 1\right\} e^{-\frac{k^2 \pi^2}{H^2} c_v \cdot t} \quad (73)$$

and, for the compression test model:

$$u(t,y) = \sum_{k=1}^{\infty} \beta_k \sin\left[\frac{(2k-1)\pi}{2H} y\right] e^{-\frac{(2k-1)^2 \pi^2}{4H^2} c_v \cdot t} \quad (74)$$

where

$$\alpha_k = a_k \frac{k\pi}{H}; \beta_k = -b_k \frac{(2k-1)\pi}{2H}; \alpha_0 = -\sum_{1}^{\infty} \alpha_k \quad (75)$$

The following two dimensionless arguments $R$ and $T$ can be derived for the compression and for the relaxation test model, respectively:

$$Y = \frac{y}{H} \text{ or } Y = \frac{y}{2H} \quad (76)$$

$$T = \frac{ct}{H^2} \text{ or } T = \frac{ct}{4H^2} \quad (77)$$

It follows from the difference of the time factors that the model law is different for the two tests. It also follows that slower dissipation time can be expected for the OCT than for the ORT in case of identical initial conditions.

**Appendix B Qualitative behaviour of the consolidation models**

Stress state in the sample
For the ORT model, the stress variation within the sample is partly rebound, partly compression as follows. The total stress $\sigma$ is equal to the sum of a steady-state and a transient component. The former - the final total stress $\sigma_\infty$ - depends on the displacement load $v_0$ and the $H$ the sample height:

$$\sigma_\infty = \frac{E_{oed} v_0}{H} \quad (78)$$

The latter is equal to the mean pore water pressure $u_{mean}$:

$$\sigma^t(t) = u_{mean}(t) = \frac{1}{H} \int_0^H u(t,y) \, dy \quad (79)$$

Effective stress decrease takes place at the sample top ($y=0$), and compression (increase in the effective stress □) takes place in the bottom ($H=0$) of the sample

$$\varepsilon^t(t,y) = -\frac{1}{E_{oed}}[u_{mean}(t) - u(t,y)] \quad (80)$$

For the OCT model, the total normal stress is constant with time. The effective normal stress increases with time, since its transient part is equal to the pore water pressure and the consolidation model predicts basically the compression.

$$\varepsilon^t(t,y) = \frac{1}{E_{oed}} u(t,y) \quad (81)$$

**Appendix C Simplification of joined model parameters**
The analytical solution of model-version $H$ contains 8 parameters, 4 out them are with non-linear dependence (Table 1). In model-version HCRT parameter $t1$ was specified, $b$ was zero and, a minor change was made in the way of the identification of the relaxation parameter s as follows. The product of $s$ and $\sigma(0)$ denoted by $sk$

$$s_k = s\,\sigma(0) \quad (82)$$

was identified instead of s. The solution depends linearly on $sk$. As a result, the number of the non-linearly dependent parameters was decreased to 2 and, the computational work was less by two orders of magnitude than the one for model H.

In model-version HCR value of $t3$, $t1$ were specified, the modified relaxation parameter ($sk$) was identified instead of s using the facts that term $\sigma(0)$ can be expressed as the linear combination of other parameters:

$$\sigma(0) = \sigma_\infty + D \quad (83)$$

where D is the mean initial pore water pressure. In model-version HCR only one among the three parameters of the relaxation part-model was identified ($sk$), value of $t3$ was set to zero. The relaxation branch of the model was left out in model-version HC (Tables 4, 5). The number of the non-linearly dependent parameters was 1 and, the computational work was less for these than for any other model-versions (Table 6).